\documentclass{amsproc}

\usepackage{amsmath}
\usepackage{amsfonts}
\usepackage{amssymb}
\usepackage{amscd}
\input{epsf.sty}
\catcode `\@=11
\def\numberbysection{\@addtoreset{equation}{section}
         \renewcommand{\theequation}{\thesection.\arabic{equation}}}
\numberbysection
\def\subsubsection{\@startsection{subsubsection}{3}%
  \normalparindent{.5\linespacing\@plus.7\linespacing}{-.5em}%
  {\normalfont\bfseries}}
\def\a{\alpha}

\def\g{\gamma}

\def\eps{\epsilon}
\def\l{\lambda}
\def\s{\sigma}
\def\t{\tau}

\def\nn{\nonumber}
\def\la{\langle}
\def\ra{\rangle}
\def\ZZ{\mathbb {Z}/2\mathbb{Z}}

\def\Sn{\mathbb{S}_n}

\def\Hom{\mathrm{Hom}}
\def\C{\mathcal{C}}
\def\F{\mathcal{F}}
\def\Frob{\mathcal{FROB}}

\def\Zz{\mathbb {Z}/2\mathbb{Z}}

\newtheorem{thm}{Theorem}[section]
\newtheorem{lm}[thm]{Lemma}
\newtheorem{prop}[thm]{Proposition}

\newtheorem{df-pr}[thm]{Definition-Proposition}
\newtheorem{df}[thm]{Definition}
\newtheorem{qu}[thm]{Question}

\begin{document}

\title[Discrete torsion, symmetric products and the Hilbert scheme]
{Discrete torsion, symmetric products and the Hilbert scheme}
\author[Ralph M. Kaufmann]{Ralph M. Kaufmann$^*$\\
University of Connecticut\\
and Max--Planck Institut f\"ur Mathematik}

\thanks{${}^*$ Partially supported by NSF grant \#0070681}
\email{ralphk@mpim-bonn.mpg.de}

\address{University of Connecticut, Department of Mathematics,
Storrs, CT, USA and
 Max--Planck Institut f\"ur Mathematik, Bonn, Germany}


\maketitle

\section*{Introduction}
Recently the understanding of the cohomology of the Hilbert scheme
of points on K3 surfaces has been greatly improved by Lehn and Sorger \cite{LS}.
Their approach uses the connection of the Hilbert scheme to
the orbifolds given by the symmetric products of these surfaces.
We introduced a general theory
replacing cohomology algebras or more generally Frobenius
algebras in a setting of global quotients by finite groups \cite{K1}. This is our
theory of group Frobenius algebras, which
are group graded non--commutative algebras whose non--commutativity is
controlled by a group action.
The action and the grading turn these algebras
into modules over the Drinfel'd double of the group ring. The appearance of
the Drinfel'd double is natural from the orbifold point of view
(see also \cite{Kir})
and can be translated into the fact that the algebra
is a $G$--graded $G$--module algebra in the following sense:
the $G$ action acts by conjugation on the grading
while the algebra structure is
compatible with the grading with respect to  left multiplication
(cf. \cite{K3,Mo}).

In the special case
of the symmetric group, we recently proved existence and
uniqueness  for the
structures of symmetric group Frobenius algebras based on a
given Frobenius algebra \cite{K2}, providing explicit formulas
for the multiplication in the algebra.

This uniqueness has to be understood up to the action of two
groups of symmetries on group Frobenius algebras called discrete torsion
and super-twisting \cite{K3}. The set of $G$--Frobenius algebras is
acted upon by both of these groups. This action only changes
some defining structures of a Frobenius algebra in a projective
manner while keeping others fixed.

Applying this result to the global orbifold cohomology of a
symmetric product, where there is a canonical choice
for the discrete torsion and super-twists, we obtain its uniqueness.

Our latest results on this topic \cite{K3} explain the origin
of these discrete degrees of freedom. In the special case of the Hilbert
scheme as a resolution of a symmetric product the
choice of sign for the metric specifies a discrete torsion cocycle
that in turn changes the multiplication by a much discussed sign.

Assembling our results which we review we obtain:

{\bf Theorem.}{\it The cohomology of $Hilb^{[n]}$, the Hilbert scheme
of n--points for a K3 surface, is the $\Sn$ invariant part of
the $\Sn$--Frobenius Algebra
associated to the symmetric product of the cohomology of the surface
twisted by a discrete torsion. Or
in other words the unique $\Sn$--Frobenius Algebra structure for the extended
global orbifold cohomology twisted by the specific discrete torsion
which is uniquely determined by the map of \cite{LS}.
In general, the sequence of spaces $Hilb^{[n]}$
gives rise to the twisted second quantization of the underlying surface
on the cohomological (motivic) level.}

Here the term associated refers to the uniqueness result of
\cite{K2} stated above.

This result follows from a series of considerations which we will review.
The logic is roughly as follows:

The theoretical background for our considerations was first presented
at \cite{talk} and is given in \cite{K1}
where we showed that the algebras arising from the ``stringy'' study
of objects with a global action by a finite group are so--called
$G$--Frobenius algebras. These algebras are non--commutative group
graded extensions of their classical
counterpart, Frobenius algebras, which arise for instance
in the study of manifolds
as cohomology algebras and in the study of singularities with an isolated
critical point as a Milnor ring. A $G$--action by automorphisms
is part of the data of a $G$--Frobenius algebra and taking invariants
under this action yields a commutative algebra.

Given an object, such as a manifold, together with a finite group action
and a functor, such as cohomology, one would like to augment this
functor to take values in $G$--Frobenius algebras.
The underlying additive structure of the $G$--Frobenius algebra
is given by evaluating the functor on each fixed point set for
each group element and forming the direct sum.
This yields a collection of Frobenius algebras,
one for each group element. The Frobenius algebra for the identity
element is the Frobenius algebra associated to the object itself and
is called the identity sector. For the other algebras, called twisted sectors,
we only retain their structure as modules over the identity sector, together
with their pairings -- the module structure over the identity sector
being induced by inclusion maps. Furthermore, there is a $G$ action on the
identity sector.

Any ``stringy'' extension of the original functor will respect
these structures, and add a {\em group graded multiplication} and
a group action by automorphisms on the whole algebra which is
compatible with the group action on the identity sector.

It is possible to classify all such ``stringy'' extensions
in the special case when all the twisted sectors are cyclic modules over
the identity sector. Such $G$--Frobenius algebras are called special.
 The classification is in terms of group cohomological
data as shown in \cite{K1}. The cyclicity condition is met in the situation
of singularities with an isolated critical point at zero as well
as for symmetric products which are the global quotient of the
$n$-th power by the $n$th symmetric group $\Sn$. For the $n$th symmetric
product of an object the untwisted sector is the $n$th tensor power of the
Frobenius algebra of this object, while the twisted sector for
a permutation is again a tensor power of the Frobenius algebra of this
object, but to the power of the number of cycles in the permutation.
It can be checked that these twisted sectors are indeed cyclic.

Imposing the cyclicity condition (i.e.\ restricting to special
$\Sn$ Frobenius algebra) and an additional grading condition it is
possible to make the classification concrete. The additional
grading condition is satisfied in the case of symmetric products.

First, we showed in \cite{K2}, if such a structure exists it is essentially
unique.

This uniqueness is essential in the following sense: as explained
in \cite{K3} given one ``stringy'' extension of the data
considered above, it is possible to produce another extension with
the use of a group cocycle in $Z^2(G,k^*)$. In the setting of
super ($\Zz$--graded) algebras, we can also produce yet another
extension for each element of ${\rm Hom}(G,\Zz)$. These twists of
the original extension can be achieved via a tensor product with a
twisted group ring or a super--graded group ring \cite{K3}. This
yields an action of both groups $Z^2(G,k^*)$ and ${\rm
Hom}(G,\Zz)$. These actions are called discrete torsion and super
twist, respectively. Thus essentially unique means unique up to
the action of these two groups.

Second, in the case of symmetric products, the unique structure
exists, as we showed in \cite{K2}. There is also a canonical
choice of an initial algebra structure upon which discrete torsion
and super--twists act.

The proof of existence relies on a general formalism which makes
use of the fact that in a setup such as a manifold with a finite
group action, we can also regard the fixed point sets of all
subgroups generated by several elements of the group. These fixed
point sets are then the intersection of the fixed point sets for
the individual generators. The general setting in which this is
possible is the setting of intersection $G$--Frobenius algebras.
In this framework, one can show that the multiplication factors
through double intersections while the associativity equation is
to be checked on triple intersections.

In the case where one is considering the symmetric product of
a manifold, the canonical non--commutative structure coincides with
the one found by \cite{FG} and its commutative
invariants are those of the Chen--Ruan orbifold
cohomology as calculated in \cite{U}. Notice that our uniqueness result
makes no reference to any space of maps or to any specific ``stringy'' extension,
 but only depends on the algebraic structure and is thus common to all ``stringy''
extensions.

Another non--commutative ``stringy'' multiplication based on
the additive data underlying symmetric products was given in \cite{LS}
in their consideration of Hilbert schemes of K3 surfaces. By our result
this has to be related to the one stemming from symmetric products by
either a twist by discrete torsion or a super--twist.
Indeed there is a twist by discrete torsion, which produces the
algebra of \cite{LS}. This discrete torsion cocycle is actually trivial
on the level of cohomology, but naturally induces a sign change for the
multiplication and the metric which is  given in \cite{LS} on the level
of group invariants. From our considerations of the action of
discrete torsion \cite{K3} this cocycle is actually already fixed
by the choice of sign for the metric, which by geometric reasoning
(resolution of $A_2$ singularities) has to be negative definite.

Lastly, the family of multiplications found in \cite{QW}
can be identified as the family of arising from
twisting with discrete torsion cocycles which
preserve the grading condition.

The paper is organized as follows: In \S1 we present the
general functorial setup for extending functors to
Frobenius algebras to those with values in
$G$--Frobenius algebras. \S2 contains the basic definitions
of $G$--Frobenius algebras and special $G$--Frobenius algebras \cite{K1},
for which the possible extensions are classifiable in terms
of group cohomological data.
In \S3 we introduce intersection Frobenius algebras which
are adapted to the situation in which one can take successive intersections
of fixed point sets. This structure is needed in order to show the
existence of symmetric products and in general it is shown how in
such a situation the multiplication can be defined via double intersections
while the associativity equations are naturally given by total symmetry
in the triple intersection.
\S4 reviews our analysis of discrete torsion \cite{K3}.
\S5 recalls our results on the structure of $\Sn$--Frobenius algebras \cite{K2}
and contains the existence and  uniqueness statements.
\S6 assembles these results in the  case of any $\Sn$--Frobenius algebra
twisted by a specific discrete torsion uniquely determined from the map \cite{LS}.
This result applied to
the situation of the Hilbert scheme yields the theorem above.

\section*{Notation} For the remainder of the paper $G$ is a fixed
finite group and $k$ is a field of characteristic $0$.

\section*{Acknowledgments}
I would like to thank the Max--Planck--Institut for its kind
hospitality and also gratefully  acknowledge the support from the
NSF. I would also like to thank the organizers of the conference.
Special thanks go to Takashi Kimura, since it was a discussion
with him that was the initial spark for our approach to discrete
torsion. Most importantly I wish to thank Yuri I.\ Manin whose
deep insights into the beauty of mathematics  have been
a continuous source of inspiration.

\section{Functorial setup}

\subsection{General background}
We will consider objects  $X$ together with the
 action of a finite group $G$. In this situation, one would classically
study the invariants or the quotient of $X$ by $G$. In stringy
geometry for global quotients, however, it is the aim to enlarge
this picture to consider the fixed point sets for all group
elements, together with an induced $G$ action on them. Here $G$
acts on the fixed point sets by conjugation of the group elements
labelling the fixed point sets.  The classical part is then
represented by $X$ considered as the fixed point set of the
identity in $G$ and the $G$ action on this fixed point set.

In particular, classical functors such as cohomology which takes values
in Frobenius algebras should have an augmented counterpart including
the information about all the fixed point sets. The augmented functors
should take values in $G$--Frobenius algebras as defined in \cite{K1}.

Physically this can be seen as the transition from topological
field theory (TFT) to a finite
gauge group TFT (see \cite{DW,DVVV,K1}).

The functorial setup of extending from functors with values in Frobenius
algebras to those with values $G$  Frobenius algebras in the following:

Let $\Frob$ be the category of Frobenius algebras, whose objects are Frobenius
algebras and morphisms are maps which respect all the structures.

\begin{df}{A  $G$--category}
is a category $\C$ where for each object
$X \in Ob(\C)$ and each $g \in G$ there exists an object $X^g$ and a
 morphism
$i_g \in \Hom(X^g,X)$  with $X^e=X$ and $i_e= id$
and there are isomorphisms $\psi_{g,g^{-1}}\in \Hom(X^{g},X^{g^{-1}})$.

We call a category a $G$ intersection category
if it is a $G$ category and for each pair
$(g,h) \in G\times G$ and object $X\in Ob(\C)$ there are isomorphisms
$\psi\in \Hom( (X^g)^h, (X^h)^g)$ and morphisms
$i^{gh}_{g,h} \in \Hom((X^g)^h,X^{gh})$.

A {$G$--action} for a $G$--category is given by
a collection of morphisms $\phi_g(X,h) \in \Hom(X^h,X^{ghg^{-1}})$
which are compatible with the structural morphisms and
satisfy $\phi_g(X,g'hg^{\prime -1})\phi_{g'}(X,h)= \phi_{gg'}(X,h)$.
\end{df}

\subsection{Examples}
Examples of an intersection $G$--category with $G$--action are
categories of spaces equipped  with a $G$--action whose fixed
point sets are in the same category. Actually this is the category
of pairs $(X,Y)$ with $X$ say a smooth space with a $G$--action
and $Y$ a subspace of $X$. Then $(X,Y)^g:=(X,Y\cap Fix(g,X))$ with
$Fix(g),X$ denoting the fixed points of $g\in G$ in X, and
$i_g=(id,\iota_g)$ with $\iota_g: Y\cap Fix(g)\rightarrow Y)$
being the inclusion. It is enough to consider pairs $(X,Y)$ where
$Y\subset X$ is the set fixed by a subgroup generated by an
arbitrary number of elements of $G$: $H:=\langle
g_1,\dots,g_k\rangle$

We could also consider the action on the $X^g$ to be trivial and
set $(X^g)^h:= X^g$. This will yield a $G$--category.

Also the category of functions $f:{\bf C}^n\rightarrow {\bf C}$
with an isolated singularity at $0$ together with a group action
of $G$ on the variables induced by a linear action of $G$ on the
linear space fixing the function is an example of a $G$--category.
This is a category of triples $({\bf C}^n, f: {\bf C}^n
\rightarrow {\bf C}), \rho \in {\rm Hom}(G,GL(n))$ such that $f$
has an isolated singularity at zero and $f(\rho({\bf z}))= f({\bf
z})$ for ${\bf z}\in {\bf C}^n$ with morphisms being linear
between the linear spaces such that all structures are compatible.
The functor under consideration is the local ring or Milnor ring.
Again we set $(X^g)^h:= X^g$.

 Here the role of
the fixed point set is played by the linear fixed point set and
the restriction of the function to this fixed point set
(cf.\cite{K1}). Again we can consider pairs of and object and a
subobject as above in order to get an intersection $G$--category.

Our main examples are smaller categories such as a global
orbifold. As a $G$ category, the objects are the fixed point sets
of the various cyclic groups generated by the element of $G$ and
the morphisms being the inclusion maps. Again we set $(X^g)^h:=
X^g$. For a global orbifold, we can also consider all fixed point
sets of the groups generated by any number of elements of $G$ as
objects together with the inclusion maps as morphisms. This latter
will render a $G$--intersection category.

The same is true for isolated singularities. Here the objects are
the restriction of the function to the various subspaces fixed by
the elements of $g$ together with inclusion or for the $G$--intersection
category we consider all intersections of these subspaces together with the
restriction of the function to these subspaces as objects,
again with the inclusion morphisms.

\subsection{The classification/reconstruction program}
Given a functor to Frobenius algebras (like cohomology),
we would like to find its stringy counterpart for global quotients.

Now, suppose we have a
$G$--category $\C$ and a contravariant functor $\mathcal{F}$
from $\C$ to $\Frob$. In this setting
there might be several schemes to define a ``stingy geometry'' by
augmenting the functor to take values in $G$--Frobenius algebras.
But all of these schemes have to have the same additive structure provided
by the ``classical orbifold picture'' (see \ref{classical}) and satisfy
the axioms of $G$--Frobenius algebras (see \S 2). Furthermore there
are more structures which are already fixed in this situation,
which is explained below. These data can sometimes be used to classify
the possible algebra structures and reconstruct it when the classification
data is known. In the case of so--called special $G$--Frobenius algebras
a classification in terms of group cohomology classes is possible.

There are some intermediate steps which contain partial information
that have been previously considered, like
the additive structure, dimensions etc., as discussed in \ref{classical}.

\subsubsection{The ``classical orbifold picture''}
\label{classical} Now, suppose we have a $G$--category $\C$ and a
contravariant functor $\mathcal{F}$ from $\C$ to $\Frob$, then for
each $X \in Ob(\C)$, we naturally obtain the following collection
of Frobenius algebras: $(\F(X^g):g\in G)$ together with
restriction maps $r_g = \F(i_g): \F(X) \mapsto \F(X^g)$.

One possibility is to regard the direct sum of the Frobenius algebras
$A_g:=\F(X^g)$.

The first obstacle is presented in the presence of a
grading, say by ${\bf N}, {\bf Z}$ or ${\bf Q}$; as it is well known that
the direct sum of two graded Frobenius algebras is only well defined
if their Euler dimensions  (cf.\ e.g.\ \cite{K3}) agree. This can, however,
be fixed by using the shifts
$s^+$ discussed in \ref{shifts}. If the grading was originally in ${\bf N}$
these shifts are usually in $\frac{1}{2}{\bf N}$, but in the complex case still lie in ${\bf N}$.

Furthermore, if we have a $G$--action on the $G$ category, it
will induce the structure of a $G$--module on this direct
sum.

Each of the Frobenius algebras $A_g$ comes equipped with
its own multiplication,
so there is a ``diagonal'' multiplication
for the direct sum which is the direct sum of these multiplications.

Using the shift $s^+$ it is possible to define a ``classical
theory'' by considering the diagonal algebra structure and taking
$G$--invariants. This is the approach used in  \cite{AS}, \cite{T}
and \cite{AR}. The paper \cite{AS} shows that this structure
describes the $G$--equivariant rather than the $G$--invariant
geometry.

One can of course forget
the algebra structure altogether and retain only
the additive structure. This was done e.g.\ in \cite{S} in the setting of
V--manifolds (i.e.\ orbifolds).
Concentrating only on the dimensions one arrives for instance at
the notion of ``stringy numbers'' \cite{BB}.

\subsubsection{The ``stringy orbifold picture''}

The ``diagonal'' multiplication is however {\em not} the right object to study
from the perspective of ``stringy geometry'' or a TFT with a finite
gauge group \cite{K1,CR}.
The multiplication should rather be $G$--graded, i.e.\ map
$A_g \otimes A_h \rightarrow A_{gh}$. We call such a product ``stringy'' product.

Here the natural question is the following:

\begin{qu}
Given the additive structure
of a $G$--Frobenius algebra, what are the possible ``stringy'' products?
\end{qu}

A more precise version of this question is the setting of our reconstruction program \cite{K2, K3}.

\subsubsection{The $G$--action} One part of the structure of
a $G$--Frobenius algebra is the $G$--action. If the $G$--category
is already endowed with a $G$--action, we can use it to reconstruct the
$G$--action on the $G$--Frobenius algebra, which in turn limits the
choices of ``stringy'' products to those that are compatible with it.

\subsubsection{Invariants} By definition $G$--Frobenius algebras
come with a $G$ action whose invariants form a commutative algebra.
Due to the nature of the $G$ action
this commutative algebra is graded by conjugacy classes, and under
certain conditions the metric descends and
the resulting algebra is again Frobenius. The induced multiplication
is multiplicative in the conjugacy classes and we call such
a multiplication commutative ``stringy''.

\subsubsection{Examples}
Examples of commutative ``stringy'' products
are orbifold (quantum) cohomology \cite{CR}.
For cohomology of global orbifolds
 it was shown in \cite{FG} and recently in \cite{JKK} that
there is a group graded version for global orbifold cohomology
which has the structure of a $G$ Frobenius algebra, as we had previously
postulated \cite{talk}.
For new developments on quantum
deformations of the $G$--Frobenius algebras see \cite{JKK}.

\subsubsection{Special $G$--Frobenius algebras}
The special reconstruction data reflects this situation in the case
that the $A_g$ algebras are cyclic $A_e$ modules. This is a restriction
which leads to an answer in terms of cocycles for a large
class of examples. This class includes all Jacobian Frobenius algebras as
well as symmetric products and special cases of geometric actions on
manifolds (were the cohomology of the fixed point sets is generated by restriction
from the ambient space).

The general idea can be generalized to the
non--cyclic case, although computations get more involved.

\begin{df}
{Given a $G$--category $\C$, we call the tuple $(X^g):g\in G$ a G--collection.

The category of $G$--collections of a $G$--category is the category whose
objects are $G$--collections and whose morphisms are collections of morphisms
$(f^g)$
s.t.\ the diagrams
$$\begin{matrix}
X^g&\stackrel{i_g}{\rightarrow}&X\\
\downarrow f^g&&\downarrow f\\
Y^g& \stackrel{i_g}{\rightarrow} & Y
\end{matrix}$$
commute.}
\end{df}

\begin{df}
{ A G--Frobenius functor is a functor from the category of $G$--collections
of a $G$--category to $G$--Frobenius algebras.}
\end{df}

\subsection{Reconstruction/classification}
The main question of the reconstruction/classification
 program is whether one can extend
a functor from a $G$--category $\C$ to Frobenius algebras to a
$G$--Frobenius functor, and if so how many ways are there to do this.

One can view this as the analogue of solving the associativity equations for
general Frobenius algebras. Some of the solutions correspond to
quantum cohomology, some to singularities, etc. and maybe others
to other ``string''--schemes. The structures of possible
``stringy'' products provide a common approach. The systematic
consideration of all possible products confines the choices of
string equivalents of classical concepts and allows to identify divers
approaches.

The answer to the main
question of reconstruction/classification can be
answered in the special case where all of the twisted
sectors are cyclic in terms of group cohomological data (see below).
This is the content of the Reconstruction Theorem of \cite{K1}.

The consequences are sometimes quite striking as in the case
of symmetric products, where there is only {one}
possible ``stringy'' orbifold
product.

The restrictions on the possible multiplicative structures are even
stricter if one is considering data stemming from a $G$--intersection category.

\section{$G$--Frobenius algebras}
\label{orb}

We fix a finite group $G$ and denote its unit element by $e$. We
furthermore fix a ground field $k$ of characteristic zero for
simplicity. With the usual precautions the characteristic
of the field does not play an important role and furthermore
the group really only needs to be completely disconnected.

\begin{df}
{ A {G--twisted Frobenius algebra} ---or $G$--Frobenius algebra for short---
over
a field  $k$ of characteristic 0 is
$<G,A,\circ,1,\eta,\varphi,\chi>$, where

\begin{tabular}{ll}
$G$&finite group\\
$A$&finite dim $G$-graded $k$--vector space \\
&$A=\oplus_{g \in G}A_{g}$\\
&$A_{e}$ is called the untwisted sector and \\
&the $A_{g}$ for $g \neq
e$ are called the twisted sectors.\\
$\circ$&a multiplication on $A$ which respects the grading:\\
&$\circ:A_g \otimes A_h \rightarrow A_{gh}$\\
$1$&a fixed element in $A_{e}$--the unit\\
$\eta$&non-degenerate bilinear form\\
&which respects grading i.e. $\eta|_{A_{g}\otimes A_{h}}=0$ unless
$gh=e$.\\
\end{tabular}

\begin{tabular}{ll}
$\varphi$&an action  of $G$ on $A$
(which will be  by algebra automorphisms), \\
&$\varphi\in \mathrm{Hom}(G,\mathrm{Aut}(A))$, s.t.\
$\varphi_{g}(A_{h})\subset A_{ghg^{-1}}$\\
$\chi$&a character $\chi \in \mathrm {Hom}(G,k^{*})$ \\

\end{tabular}

\vskip 0.3cm

\noindent Satisfying the following axioms:

\noindent{\sc Notation:} We use a subscript on an element of $A$ to signify that it has homogeneous group
degree  --e.g.\ $a_g$ means $a_g \in A_g$-- and we write $\varphi_{g}:= \varphi(g)$ and $\chi_{g}:= \chi(g)$.

\begin{itemize}

\item[a)] {Associativity}

$(a_{g}\circ a_{h}) \circ a_{k} =a_{g}\circ (a_{h} \circ a_{k})$
\item[b)] {Twisted commutativity}

$a_{g}\circ a_{h} = \varphi_{g}(a_{h})\circ a_{g}$
\item[c)]
{$G$ Invariant Unit}:

$1 \circ a_{g} = a_{g}\circ 1 = a_g$

and

$\varphi_g(1)=1$
\item[d)]
{Invariance of the metric}:

$\eta(a_{g},a_{h}\circ a_{k}) = \eta(a_{g}\circ a_{h},a_{k})$

\item[i)]
{Projective self--invariance of the twisted sectors}

$\varphi_{g}|A_{g}=\chi_{g}^{-1}id$

\item[ii)]
{$G$--Invariance of the multiplication}

$\varphi_{k}(a_{g}\circ a_{h}) = \varphi_{k}(a_{g})\circ  \varphi_{k}(a_{h})$

\item[iii)]

{Projective $G$--invariance of the metric}

$\varphi_{g}^{*}(\eta) = \chi_{g}^{-2}\eta$

\item[iv)]
{Projective trace axiom}

$\forall c \in A_{[g,h]}$ and $l_c$ left multiplication by $c$:

$\chi_{h}\mathrm {Tr} (l_c  \varphi_{h}|_{A_{g}})=
\chi_{g^{-1}}\mathrm  {Tr}(  \varphi_{g^{-1}} l_c|_{A_{h}})$
\end{itemize}}
\end{df}

\subsection{Remark} In the case of trivial characters the notion
of $G$--Frobenius algebras has appeared also under the name
of group--crossed algebras in \cite{Tu}, where they appeared
from the point of view of homotopy field theory.

\subsection{$G$--graded tensor product} Given two $G$--Frobenius algebras
$$\la G,A,\circ,1,\eta,\varphi,\chi\ra \ \ \text{ and }
\la G,A',\circ',1',\eta',\varphi',\chi'\ra$$ we
defined \cite{K1} their tensor product
as $G$--Frobenius algebras to be the $G$--Frobenius algebra

$\la G,\bigoplus_{g \in G}( A_g \otimes A'_g),
\circ\otimes \circ',1\otimes 1',\eta\otimes \eta',\varphi\otimes \varphi',
\chi\otimes \chi'\ra$.

We will use the short hand notation $A \hat \otimes A'$ for this product.

For the program outlined in \S1 the following data is the starting point
in order to construct a $G$--Frobenius algebra.
\begin{df}
A reconstruction data is collection of Frobenius algebras
$$(A_g,\eta_g,1_g):g\in G$$ together with  maps of algebras $r_g:
A_e \rightarrow A_g$, isomorphisms $\psi_g:A_g \tilde\rightarrow
A_{g^{-1}}$ and a $G$--action $\varphi$ on $A_e$.
\end{df}

In general, we would like to find the $G$--Frobenius algebra structures
compatible with these data. For many purposes such as symmetric
products it is however enough to restrict to a more specialized situation.

\subsection{Special $G$--Frobenius algebras}
We briefly review special $G$--Frobenius algebras. For this class
of algebras which include the algebras having their origin in singularities
with isolated singularities and symmetric products a classification
of all possible stringy multiplications is possible in terms of
group cohomological data. For details
see \cite{K1,K2}.

\begin{df}
We call a $G$-Frobenius algebra special if all $A_g$ are cyclic
$A_e$ modules via the multiplication $A_e \otimes A_g \rightarrow A_g$
and there exists a collection of cyclic generators $1_g$ of $A_g$ such that
$\varphi_g(1_h)= \varphi_{g,h}1_{ghg^{-1}}$ with $\varphi_{g,h}\in k^*$.
\end{df}

The last condition is automatic, if the Frobenius algebra $A_e$
only has $k^*$ as invertibles, as is the case for cohomology algebras of
connected compact manifolds
and Milnor rings of quasi--homogeneous functions with an
isolated critical point at zero.

Fixing the generators $1_g$ we obtain maps $r_g:A_e \rightarrow A_g$ by setting
$r_g(a_e)= a_e1_g$. This yields a short exact sequence

\begin{equation}
0\rightarrow I_g \rightarrow A_e \stackrel{r_g}{\rightarrow} A_g \rightarrow 0
\end{equation}

It is furthermore useful to fix a section $i_g$ of $r_g$.

We denote the concatenation $\pi_g:= i_g \circ r_g$.

\begin{df}
A special $G$ reconstruction datum is a collection of Frobenius
algebras $(A_g,\eta_g,1_g): g\in G$ together with
an action of $G$ by algebra automorphisms
on $A_e$ and the structure of a cyclic $A_e$ module algebra on each $A_g$ with
generator $1_g$ such that $A_g$ and
$A_g^{-1}$ are isomorphic as of $A_e$ modules algebras.
\end{df}

\begin{df}
{\it Given a Frobenius algebra $A_e$ and a collection of
cyclic $A_e$--modules $A_g:g \in G$
a graded cocycle  is a map $\g: G\times G \rightarrow A_e$
which satisfies
$$\g(g,h)\g(gh,k)\equiv \g(g,hk)\g(h,k) \; \mathrm{ mod }\;  I_{ghk}$$

Such a cocycle is called
section independent if
$$(I_g + I_h)\g(g,h) \subset I_{gh}$$

Two such cocycles are considered to be the same if $\g(g,h) \equiv
\g'(g,h) \; \mathrm{ mod }\; I_{gh}$ and isomorphic, if they are
related by the usual scaling for group cocycles.

Given non--degenerate parings $\eta_g$ on the $A_g$,
a cocycle is said to be compatible with the metric, if
$$
\check r_g(1_g) = \g(g,g^{-1})
$$
where $\check r$ is the dual in the sense of vector spaces with
non--degenerate metric.
}
\end{df}

We will again use the notation $\g_{g,h}:=\g(g,h)$.

\subsection{Special $G$--Frobenius structure in terms of the cocycles}
\label{special}
Fixing a cyclic generator $1_g \in A_g$, a special $G$--Frobenius
algebra is completely characterized by two
cocycles:  $\g$ a section independent graded cocycle compatible with the
metric and $\varphi\in Z^1(G,k^*[G])$
where $k^*[G]$ is the group ring restricted to invertible coefficients
with $G$--module structure induced by the adjoint action:
$$
\phi(g)\cdot(\sum_h \mu_h h)= \sum_h \mu_h ghg^{-1}
$$

The multiplication and $G$--action on the generators defines these
cocycles. Set

$$
1_g 1_h = \g_{g,h} 1_{gh} \quad \varphi_{g}(1_{h}) = \varphi_{g,h}
1_{ghg^{-1}}
$$

Defining $\varphi(g) := \sum_h \varphi_{g,h} ghg^{-1}$ and
 $\g(g,h):=\g_{g,h}$ we obtain the desired cocycles.

The section independence follows from the fact that

$$(I_g+I_h)\g_{g,h}1_{gh}= (I_g+I_h)1_g1_h=0$$

In general, the multiplication is thus given by
\begin{equation}
\label{specialmult}
a_g b_h = i_g(a_g)i_h(b_h)\g_{g,h}1_{gh}
\end{equation}
for any choice of sections $i_g$.

The cocycles furthermore satisfy the following two compatibility equations:
\begin{equation}
\label{grpcompat}
\varphi_{g,h}\g_{ghg^{-1},g} = \g_{g,h}
\end{equation}
and
\begin{equation}
\label{algaut}
\varphi_{k,g} \varphi_{k,h} \g_{kgk^{-1},khk^{-1}}
= \varphi_{k} (\g_{g,h}) \varphi_{k,gh}
\end{equation}
which follow from the twisted commutativity and
the fact that $\varphi$ acts by automorphisms.

\subsubsection{Remark}
Notice that if $\g_{g,h}$ is non--zero i.e. $A_gA_h \neq 0$ then
(\ref{grpcompat})
determines $\varphi_{g,h}$. We also would like to remark that
 if $A_gA_hA_k \neq 0$
(\ref{algaut}) follows from  (\ref{grpcompat})
(cf.\ \cite{K1}).

\begin{df}
We call a pair of a section independent cocycle and
a non--abelian cocycle compatible if they satisfy the equations
(\ref{grpcompat}) and (\ref{algaut}).
\end{df}
\begin{thm}(Reconstruction) Given a special $G$ reconstruction datum
the structures of special $G$--Frobenius algebras are in 1--1
correspondence compatible pairs of
a graded, section independent $G$ 2--cocycle with values in $A_e$ that is
compatible with the metric and a
non--abelian $G$ 2--cocycle with values in $K^{*}$.
Satisfying the following conditions:
\begin{itemize}
\item[i)]$\varphi_{g,g}=\chi_g^{-1}$
\item[ii)]
$\eta_{e}(\varphi_{g}(a),\varphi_{g}(b)) =
    \chi_{g}^{-2}\eta_{e}(a,b)$
\item[iii)] The projective trace axiom
$\forall c \in A_{[g,h]}$ and $l_c$ left multiplication by $c$:
\begin{equation}
\chi_{h}\mathrm {Tr} (l_c  \varphi_{h}|_{A_{g}})=
\chi_{g^{-1}}\mathrm  {Tr}(  \varphi_{g^{-1}} l_c|_{A_{h}})
\end{equation}
\end{itemize}
\end{thm}

\subsection{Remark}
Changing the cyclic generators by elements
of $k^*$ leads to isomorphic $G$--Frobenius algebras and
to cohomologous cocycles $\g,\varphi$ in $Z^2(G,A_e)$ and $Z^2(G,k^*[G])$.

\subsection{Grading and Shifts}
\label{shifts}
Consider the graded version for Frobenius algebras,
where each Frobenius $A_g$ algebra comes naturally graded,
 e.g.\ by cohomological degree or weight in the case
of a quasi--homogeneous isolated singularity. Usually this grading
takes values in ${\bf Q}$.
In this case the metric is also of a fixed degree, e.g.\ the dimension
or the highest weight in the Milnor ring.

Then each $A_g$ is graded as well. We denote the degree of the
pairing $\eta_g$ by $d_g$ and  also use the shorthand notation
$d:= d_e$.

\subsubsection{Remark}
It is well known that the direct sum of graded Frobenius algebras is a graded
Frobenius algebra only if the degrees match (see e.g. \cite{K1}).

\subsubsection{The shifts}
From the previous Remark it follows that in order to form the sum
$A:=\bigoplus_{g\in G}A_g$, we need to shift the degrees of the
elements of $A_g$ at least uniformly, i.e.\ if an element $a_g$ in
$A_g$ has degree $\deg(a_g)$ we assign to it the new shifted
degree $\deg^s(a_g) = \deg(a_g) + s_g$.

 This observation does not
fix the shifts uniquely. Let us denote the shift in degree of
$A_g$ by $s_g$ and set

$$s_{g}^+ := s_{g}+s_{g^{-1}}, \quad s_{g}^-:= s_{g}-s_{g^{-1}}$$

Then
$$s_{g}^+:= d-d_{g}$$
for the grading reasons mentioned above.

The shift $s^-$ is not fixed, however, there is a standard choice
provided there exists a canonical choice of linear representation of $G$.

\begin{df}
{\it The standard shift for a G--Frobenius algebra with a choice of linear
representation $\rho: G \rightarrow GL_n(k)$
is given by
$$s_{g}^+:= d-d_{g}$$
and
\begin{multline*}
s_{g}^- := \frac{1}{2\pi i}\mathrm{tr} (\log(g))-\mathrm{tr}(\log(g^{-1})):=
\frac{1}{2\pi i}(\sum_i \l_i(g)-\sum_i \l_i(g^{-1}))\\
=\sum_{i: \l_i \neq 0} (\frac{1}{2\pi i}2\l_i(g)-1)
\end{multline*}
where the $\l_i(g)$
are the logarithms of the eigenvalues
of $\rho(g)$ using the branch with arguments in $[0,2\pi)$ i.e.\
cut along the positive real axis.}
\end{df}

In total we obtain:

$$
s_{g}= \frac{1}{2}(s_g^+ + s_g^-)= \frac{1}{2}(d-d_g)
+ \sum_{i:\l_i \neq0} (\frac{1}{2\pi i}\l_i(g)-\frac{1}{2})
$$

\subsubsection{Remark} This grading having its origin in physics
specializes to the so--called age grading or the orbifold grading of \cite{CR}
in the respective situations.

\subsection{Super-grading}
\label{super} We can enlarge the framework by considering
super--algebras rather than algebras. This will introduce the
standard signs.

The action of $G$ as well as the untwisted sector should be even.
The axioms that change are

\begin{itemize}

\item[b$^{\sigma}$)] {Twisted super--commutativity}

$a_{g}\circ a_{h} = (-1)^{\tilde a_g\tilde a_h} \varphi_{g}(a_{h})\circ a_{g}$

\item[iv$^{\sigma}$)]
{Projective super--trace axiom}

$\forall c \in A_{[g,h]}$ and $l_c$ left multiplication by $c$:

$\chi_{h}\mathrm {STr} (l_c  \varphi_{h}|_{A_{g}})=
\chi_{g^{-1}}\mathrm  {STr}(  \varphi_{g^{-1}} l_c|_{A_{h}})$
\end{itemize}
where $\mathrm{STr}$ is the super--trace.

Here we denoted by $\tilde a$ the $\mathbb{Z}/2\mathbb{Z}$ degree of $a$.

\section{Intersection G--Frobenius algebras}

We will now concentrate on the situation of functors from $G$--intersection
categories to Frobenius algebras.

Given a $G$--class in such a category a functor to Frobenius algebras will
provide the following structure which reflects the possibility
to take fixed point sets iteratively. Say we look at the fixed
points with respect to
elements $g_1, \dots, g_n$. These fixed point sets
will be invariant under the group spanned by the elements
$g_1, \dots, g_n$ and they are just the intersection of
the respective fixed point sets of the elements $g_i$.
The underlying spaces are therefore invariant with respect
to permutation of the elements $g_i$, and if $g$ appears twice
 among the $g_i$ then one can shorten
the list by omitting one of the $g_i$. Also if a list $g_i$
includes $g^{-1}$ we may replace it by $g$. Finally, the fixed
point set under the action of the group generated by two elements
$g$ and $h$ is a subset of the fixed point set of the group
generated by their product $gh$. Translating this into the
categorical framework, we obtain:

\begin{df}
A $G$--intersection Frobenius datum of level $k$ is the following:
For each collection $(g_1,\ldots, g_n)$ with $n\leq k$ of elements
of $G$, a Frobenius algebra $A_{g_1,\dots,g_n}$ and the following
maps:

Isomorphisms
$$\Psi_{\sigma}:
 A_{g_1,\dots,g_n}\rightarrow A_{g_{\sigma(1)},\dots,g_{\sigma(n)}}$$
for each $\sigma \in \Sn$ called {permutations}.

Isomorphisms
$$\Psi^{g_1,\dots, g_i, \dots, g_n}_{g_1,\dots, g_i^{-1},\dots ,g_n  }:
A_{g_1,\dots, g_i, \dots, g_n} \rightarrow
 A_{g_1,\dots, g_i^{-1},\dots ,g_n}$$
 commuting with the permutations.

Morphisms
$$r_{g_1,\dots, g_i, \dots g_n}^{g_1,\dots , \hat g_i,\dots ,g_n}:
A_{g_1,\dots,\hat g_i, \dots, g_{n}}
\rightarrow A_{ g_1,\dots,g_n}$$
commuting with the permutations. (Here the symbol $\hat {}$
is used to denote omission.)

Isomorphisms
$$i_{g_1,\dots, g, \dots , g,\dots, g_n}^{g_1,\dots g, \dots, \hat g,\dots ,g_n}:
A_{g_1,\dots, g, \dots , g,\dots, g_n}
\rightarrow A_{g_1,\dots g, \dots, \hat g,\dots ,g_n}$$
commuting with the permutations.

And finally morphisms:
$$r^{g_1,\dots,g_{i}g_{i+1},\dots, g_n}_{g_1,\dots,g_{i},g_{i+1},\dots, g_n}:
A_{g_1,\dots,g_{i}g_{i+1},\dots, g_n} \rightarrow
A_{g_1,\dots,g_{i},g_{i+1},\dots, g_n}$$
commuting with the permutations.

If this data exists for all $k$ we call the data simply
$G$--intersection Frobenius datum.
\end{df}
\subsection{Notation}
We set $r_{g_1,\dots,g_n}:= r_{g_1, \dots, g_n}^{g_1,\dots, g_{n-1}}
\circ \dots \circ r_{g_1}$
and we set $I_{g_1,\dots,g_n}:= \mathrm{Ker}( r_{g_1,\dots, g_n})$.
Notice that this definition of $I_{g_1,\dots,g_n}$ is independent
of the order of the $g_i$.

\begin{df}
An intersection $G$--Frobenius algebra of level $k\geq 2$ is
an intersection $G$--Frobenius datum of level $k\geq 2$
together with a $G$--Frobenius algebra structure on $A:= \bigoplus A_g$.

An intersection $G$--Frobenius algebra of level $k\geq 2$ is called special,
if all of the $A_{g_1,\dots g_n}$ are cyclic $A_e$ module algebras generated by
the $1_{g_1, \dots g_n}$.
\end{df}

\subsection{Remarks}
\begin{itemize}
\item[1)]
In order to (re)--construct
a suitable multiplication on $\bigoplus A_g$ it is often
convenient to use the double and triple intersections (i.e.\  level 3).
Where the double intersection are used for the multiplication and
triple intersections are used to
show associativity.
\item[2)] We can use the double intersections to define $G$--Frobenius
algebras based on each of the $A_g$ i.e.\ on
 $\bigoplus_{h\in Z(g)} A_{g,h}$ for each fixed
$g$--where $Z(g)$ denotes the centralizer of $g$.
\end{itemize}

\begin{df}
A {$G$--action} for an intersection $G$--Frobenius datum of
level $k$ is given by
a collection of morphisms
$$\phi_g(A_{g_1,\dots,g_n},h) \in
\Hom(A_{g_1,\dots,g_n,h},A_{g_1,\dots,g_n,ghg^{-1}})$$
which are compatible with the structural homomorphisms and
satisfy
$$\phi_g(A_{g_1,\dots,g_n},g'hg^{\prime -1})\phi_{g'}
(A_{g_1,\dots,g_n},h)= \phi_{gg'}(A_{g_1,\dots,g_n},h)$$
\end{df}

\begin{df}
We call an intersection $G$ Frobenius datum
a special $G$ intersection Frobenius datum
datum, if all of the $A_{g_1,\dots,g_n}$ are cyclic $A_e$
module algebras via the restriction maps such that the $A_e$ module
structures are compatible with the restriction morphisms $r$. Here
the generators are given by $r_{g_1,\dots, g_n}(1)$ and the $A_e$
module structure is given by $a\cdot b:= r_{g_1,\dots,g_n}(a)b$.
\end{df}

\subsection{Remark}
In the case of special intersection $G$--Frobenius algebras, there
are two ways to look at the multiplication. One way is to use the
restrictions $r_g$ and sections $i_g$ to define the multiplication
as discussed in \S \ref{special}. A second possibility is to use
the intersection structure. This can be done in the following way:
first push forward to double intersections, second use the
Frobenius algebra structure there to multiply, then pull the
result back up to the invariants of the product, but allowing to
multiply with an obstruction class before pulling back. This is
discussed below in \S \ref{multiplication}.

The precise relation between the two procedures is given by
the following Proposition and \ref{specialmult}.

\begin{prop}\cite{K2}
{\it \label{intersect} Given a special $G$ intersection algebra datum
(of level $2$),
the following decomposition holds for section independent cocycles
$\g$:
\begin{equation}
r_{gh}(\g_{g,h}) =\check r_{g,h}^{gh}(\tilde\g_{g,h})
= i_{g,h}^{gh}(\tilde \g_{g,h})  \check r_{g,h}^{gh}(1_{g,h})
=\bar \g_{g,h}  \g_{g,h}^{\perp}
\end{equation}
 for some section $i_{g,h}$ of $r_{g,h}$, $\tilde \g_{g,h} \in (A_{g,h})^{e}$,
$\bar \g_{g,h}\in i_g,h)(A_g,h)$ of degree $e$.
and $ \g_{g,h}^{\perp}:=\check r_{g,h}^{gh}(1_{g,h})$
Here $e=s_{g}+s_{h}-s_{gh}-s^{+}_{g,h}+s^{+}_{gh}$ with $s^{+}_{g,h} :=
d-d_{g,h}$ and $d_{g,h}=\deg(\rho_{g,h})$
and we again used the unshifted degrees. (In particular if  the
$s^{-}=0$ then $e= \frac{1}{2}(s_g^++s_h^++s_{gh}^+)-s_{g,h}^+
=\frac{1}{2}(d-d_g-d_h-d_{gh})+d_{g,h}$)}
\end{prop}

Here $\check r$ is again the dual for maps between
vector spaces with non--degenerate
bilinear forms.

That is  using the multiplication in $A_{g,h}$

\begin{equation}
a_g \circ b_h =
\check r^{gh}_{g,h}(r^g_{g,h}(a_g)r^h_{g,h}(b_h)\tilde\g_{g,h})
\end{equation}

\subsubsection{Remark} The decomposition into the terms
$\tilde \g$ and $\g^{\perp}$ can be understood as decomposing the
cocycle into a part which comes from the normal bundle of $X^{g,h}
\subset X^{gh}$ which is captured by $\g^{\perp}$ and an additional
obstruction part.

Also generalizing the fact that
\begin{equation}
I_g \g_{g} = I_g \check r_g (1_g)=0
\end{equation}
the following lemma holds:
\begin{lm}
\begin{equation}
(I_g +I_h)\g_{g,h}^{\perp}\subset I_{gh}
\end{equation}
\end{lm}

\subsection{Multiplication}
\label{multiplication} Fix a special intersection $G$ Frobenius algebra
of level at least 2. From the section independence of $\g$, we
see that the multiplication
$A_g \otimes A_h \rightarrow A_{gh}$ can be factored through
$A_{g,h}$. To be more precise, we have the following commutative diagram.

$$
\begin{matrix}
A_g\otimes A_h &\stackrel{\mu}{\rightarrow} &A_{gh}\\
\downarrow r^g_{g,h}\otimes r^h_{g,h}&&
\uparrow \check r^{g,h}_{gh}\circ l_{\tilde\g_{g,h}}\\
A_{g,h}\otimes A_{g,h}&\stackrel{\mu}{\rightarrow}&A_{g,h}
\end{matrix}
$$
where $l_{\tilde\g_{g,h}}$ is the left multiplication with $\tilde\g_{g,h}$.

Vice--versa we can use this diagram as an Ansatz for any $G$ intersection
Frobenius algebra of level at least 2.

\subsection{Associativity equations}
\label{ass}
Fix an intersection $G$ Frobenius algebra of level at least $3$ then
the associativity equations can be factored through
$A_{g,h,k}$. More precisely, we have the following commutative diagram
of restriction maps:

\begin{equation}
\label{assdiagram}
\begin{matrix}
&&&&A_{ghk}&&&&\\
&&&\swarrow&&\searrow&&&\\
A_{gh}&\rightarrow&A_{gh,k}&&\downarrow&&A_{g,hk}&\leftarrow&A_{hk}\\
\downarrow&&&\searrow&&\swarrow&&&\downarrow\\
A_{g,h}&&\rightarrow&&A_{g,h,k}&&\leftarrow&&A_{h,k}\\
\end{matrix}
\end{equation}

More technically:
Using the associativity equations for the $\g$, we set
\begin{equation}
\label{tildetripel}
r_{ghk}(\g_{g,h} \g_{gh,k}):=
\g_{g,h,k}
\end{equation}

Associativity dictates that also
\begin{equation}
r_{ghk}(\g_{h,k} \g_{g,hk})=
\g_{g,h,k}
\end{equation}

By analogous arguments as for the decomposition of the $\g_{g,h}$'s one can obtain:
\begin{equation}
\g_{g,h,k}= i_{g,h,k}^{ghk}(\tilde\g_{g,h,k}) \check r_{g,h,k}^{ghk}(1_{g,h,k})
= \check r_{g,h,k}^{ghk}(\tilde\g_{g,h,k})
\end{equation}
for some $\tilde \g_{g,h,k}\in A_{g,h,k}$.

Vice--versa having defined suitable $\tilde \g_{g,h}$
to show associativity one needs to show that
\begin{equation}
\check r_{gh,k}^{ghk}
(r_{gh,k}^{gh}(\check r^{gh}_{g,h}(\tilde \g_{g,h}))\tilde \g_{gh,k})=
  \check r_{g,h,k}^{ghk}(\tilde \g_{g,h,k})
\end{equation}
for some $\tilde \g_{g,h,k}$.
This approach is actually now independent of the setup of special $G$
intersection
Frobenius algebras, where such a decomposition is guaranteed, and
is suitable for all $G$ intersection Frobenius data respectively
intersection $G$--categories.

\section{Discrete Torsion}

\label{disc}
\subsection{The twisted group ring $k^{\a}[G]$}

Recall that given an element $\a \in Z^2(G,k^*)$
one defines the twisted group ring
$k^{\a}[G]$ to be given by the same linear structure with multiplication
given by the linear extension of

\begin{equation}
g\otimes h \mapsto \a(g,h) gh
\end{equation}
with $1$ remaining the unit element.
To avoid confusion we will denote elements of $k^{\a}[G]$ by
$\hat g$ and the multiplication with $\cdot$
Thus
$$
\hat g \cdot \hat h = \a(g,h) \widehat{gh}
$$

For $\a$ the following equations hold:
\begin{equation}
\a (g,e) = \a(e,g)=1, \qquad
\a(g,g^{-1})=\a(g^{-1},g)
\end{equation}
Furthermore
$$
\hat{g}^{-1}= \frac{1}{\a(g,g^{-1})}\widehat{g^{-1}}
$$
and
$$
\hat g\cdot \hat h\cdot\hat {g}^{-1} =
\frac{\a(g,h)\a(gh,g^{-1})}{\a({g,g^{-1})}}\widehat{ghg^{-1}}
= \frac{\a(g,h)}{\a(ghg^{-1},g)} \widehat{ghg^{-1}}=
\eps(g,h)\widehat{ghg^{-1}}
$$
with
\begin{equation}
\eps(g,h):=\frac{\a(g,h)}{\a(ghg^{-1},g)}
\end{equation}
\subsubsection{Remark}
\label{conorm} If the field $k$ is algebraically closed
we can find a representative for each class $[\a]\in H^2(G,k^*)$
which also satisfies

$$\a(g,g^{-1})=1$$

\subsubsection{The $G$--Frobenius Algebra structure of $k^{\a}[G]$}
Fix  $\a \in Z^2(G,k^*)$.
Recall from \cite{K1,K2} the following structures which turn
$k^{\a}[G]$ into a special $G$--Frobenius algebra:

\begin{eqnarray}
 \g_{g,h}=\a(g,h) &&\eta(\hat g,\widehat{g^{-1}}) =\a(g,g^{-1})\nn\\
\chi_g= (-1)^{\tilde g}
&&\varphi_{g,h}=
\frac{\a(g,h)}{\a(ghg^{-1},g)}=:\eps(g,h)
\end{eqnarray}
\subsubsection{Relations}
The $\eps(g,h)$ which are by definition given as
$$\eps(g,h):= \frac{\a(g,h)}{\a(ghg^{-1},h})$$ satisfy the equations:
\begin{eqnarray}
\eps(g,e)&=&\eps(g,g)=1\\\nn
 \eps(g_1g_2,h)&=&
\eps(g_1,g_2hg_2^{-1})\eps(g_2,h)\nn\\
 \eps(k,gh)& =& \eps(k,g)\eps(k,h)\frac{\a(kgk^{-1},khk^{-1})}{\a(g,h)}\nn\\
\eps(h,g)&=&\eps(g^{-1},ghg^{-1})
\frac{\a([g,h],h)}{\a([g,h],hgh^{-1})}
\end{eqnarray}

This yields for {commuting elements}:

\begin{eqnarray}
\label{eps}
\eps(g,e)=\eps(g,g)=1 && \eps(g,h) =\eps(h^{-1},g)=\eps(h,g)^{-1}\nn\\
\eps(g_1g_2,h)= \eps(g_1,h)\eps(g_2,h)  &&
\eps(h,g_1g_2) = \eps(h,g_1) \eps(h,g_2)
\end{eqnarray}

In the physics literature discrete torsion is sometimes defined to
be a function  $\eps$ defined on commuting elements of $G$ taking
values in $U(1)$ and satisfying the equations (\ref{eps}).

\subsubsection{Remark} It is a nice exercise to check that the
trace axiom also holds (see \cite{K1,K3}).

\subsubsection{Remark} The function $\eps$ can be interpreted as
a cocycle in $Z^1(G,k^*[G])$ where $k^*[G]$ are the elements of
$k[G]$ with invertible coefficients regarded as a $G$ module by
conjugation (cf. \cite{K1,K2}). This means in particular that on
{\em commuting elements} $\eps$ only depends on the class of the cocycle $\a$.

\subsection{The action of discrete Torsion}

\begin{df}
Given a $G$--Frobenius algebra
$A$ and an element $\a \in Z^2(G,k)$, we define the
$\a$--twist of $A$ to be the $G$--Frobenius algebra
$A^{\a}:= A \hat\otimes k^{\a}[G]$.
\end{df}

\begin{prop}
\label{defprop}
Notice that as vector spaces
\begin{equation}
\label{alphaiso}
A^{\a}_{g}= A_g \otimes k \simeq A_g
\end{equation}
Using this identification the $G$--Frobenius structures  given by
(\ref{alphaiso}) are
\begin{eqnarray}
\circ^{\a}|_{A^{\a}_{g}\otimes A^{\a}_{h}}= \a(g,h) \circ &&
\varphi^{\a}_g|_{A^{\a}_h}=\eps(g,h)\varphi_g\nn\\
\eta^{\a}|_{A^{\a}_g\otimes A^{\a}_{g^{-1}}}= \a(g,g^{-1})\eta&&
\chi_g=\chi_g
\end{eqnarray}
\subsubsection{Supergraded twisted group rings}
Fix $$\a \in Z^2(G,k^*), \s \in \Hom(G,\Zz)$$ then there is a twisted
super--version of the group ring where now the relations
read
\begin{equation}
 \hat g \hat h =  \a(g,h)\widehat {gh}
\end{equation}
and the twisted commutativity is
\begin{equation}
 \hat g \hat h = (-1)^{\s(g)\s(h)}\varphi_{g}(\hat h) \hat g
\end{equation}
and thus
\begin{equation}
 \varphi_{g}(\hat h)=
(-1)^{\s(g)\s(h)}\a(g,h)\a(gh,g^{-1}) \widehat{ghg^{-1}} =:
\varphi_{g,h} \widehat{ghg^{-1}}
\end{equation}
and thus
\begin{equation}
\eps(g,h) := \varphi_{g,h} = (-1)^{\s(g)\s(h)}\frac{\a(g,h)}{\a(ghg^{-1},g)}
\end{equation}
\end{prop}

 We would just like to remark that the
axiom iv$^{\sigma})$ of  \ref{super} shows the difference between
super twists and discrete torsion.

\begin{df}
We denote the $\a$-twisted group ring
with super--structure $\s$ by $k^{\a,\s}[G]$.
We still denote $k^{\a,0}[G]$ by $k^{\a}[G]$
where $0$ is the zero map and we denote $k^{0,\s}[G]$
just by $k^{\s}[G]$ where $0$ is the unit of the group $H^2(G,k^*)$.
\end{df}

A straightforward calculation shows
\begin{lm}
 $k^{\a,\s}[G] = k^{\a}[G]\otimes k^{\s}[G]$.
\end{lm}

\begin{lm} Let $\la G,A,\circ,1,\eta,\varphi,\chi\ra$  be
a $G$--Frobenius algebra or more generally super Frobenius algebra with
super grading $\tilde{}\;\in {\rm Hom}(A,\Zz)$
 then
$A\otimes k^{\s}[G]$ is isomorphic to the super $G$--Frobenius algebra
$\la G,A,\circ^{\sigma},1,\eta^{\sigma},\varphi^{\s},\chi^{\s}\ra$ with super grading
${}^{\sim \s}$, where
\begin{eqnarray*}
\circ^{\sigma}|_{A_g\otimes A_h}=(-1)^{\tilde g \sigma(h)}\circ
&\quad &
\varphi^{\s}_{g,h} = (-1)^{\s(g)\s(h)}\varphi_{g,h}\\
\eta_g^{\sigma}=(-1)^{\tilde g \sigma(g)}\eta_g
&\quad& \chi^{\s} = (-1)^{\s(g)}\chi_g\\
\tilde a_g^{\s}= \tilde a_g + \s(g)&&\\
\end{eqnarray*}
\end{lm}

\begin{df}
Given a $G$--Frobenius algebra
$A$ a twist for $A$ is a pair of functions
$(\l:G\times G \rightarrow k^*,\mu:G\times G \rightarrow k^*)$

Such that $A$ together with the new $G$--action
$$\varphi^{\l}(g)(a) = \oplus_h \l(g,h) \varphi(g)(a_h)$$
and the new multiplication
$$
a_g \circ^{\mu} b_h = \mu(g,h) a_g \circ b_h
$$
is again a $G$--Frobenius algebra.

A twist is called universal if it is defined for all $G$--Frobenius algebras.
\end{df}

\subsubsection{Remark} We could have started from a pair of functions
$(\l:A\times A \rightarrow k^*,\mu:G\times A \rightarrow k^*)$ in order
to projectively change the multiplication and $G$ action, but it is
clear that the universal twists (i.e.\ defined for any $G$--Frobenius
algebra) can only take into account
the $G$ degree of the elements.

\subsubsection{Remark}
These twists arise from a projectivization of the $G$--structures
induced on a module over $A$ as for instance the associated
Ramond--space (cf.\ \cite{K1}). In physics terms this means that
each twisted sector will have a projective vacuum, so that fixing
their lifts in different ways induces the twist. Mathematically
this means that the $g$ twisted sector is considered to be a Verma
module over $A_g$ based on this vacuum.

\begin{thm}
Given a (super) $G$--Frobenius algebra $A$
the universal twists are in 1--1 correspondence with
elements $\a \in Z^2(G,k^*)$  and the isomorphism classes of universal
twists are given by $H^2(G,k^*)$. Furthermore
the universal super  re--gradings are
in 1-1 correspondence with $\Hom(G,\Zz)$ and these
structures can be realized by tensoring with $k^{\s}[G]$
for $\s \in  \Hom(G,\Zz)$.
\end{thm}

Here a super re--grading is a new super grading on $A$ with which
$A$ is a super $G$--Frobenius algebra and universal means that
the operation of re--grading is defined for all $G$--Frobenius algebras.

We call the operation of forming a tensor product with $k^{\a}[G]:
\a \in Z^2(G,k^*)$ a twist by discrete torsion. The term discrete
refers to the isomorphism classes of twisted $G$--Frobenius
algebras which correspond to classes in $H^2(G,k^*)$. Furthermore,
we call  the operation of forming a tensor product with
$k^{\s}[G]:\s \in {\rm Hom}(G,\Zz)$ super--twist.

\subsection{Remark} If $k$ is algebraically closed, then in each class
of $H^2(G,k^*)$ there is a representative with $\a(g,g^{-1})=1$.
Using these representatives it is possible to twist a special
$G$--Frobenius algebra without changing its underlying special
reconstruction data.

\section{Symmetric group Frobenius algebra}
In this section, we will consider the structure of special
$G$--Frobenius algebras when the group is a symmetric group $\Sn$.
The symmetric groups have two characteristics, which we will use.
First they are generated by self inverse transpositions and second
there is a natural grading on the elements given by the minimal
number of transpositions needed to present an element.

Using the latter condition in order to fix a grading compatibility for
special $\Sn$ Frobenius algebras, we showed \cite{K2} that the structure
of a special  $\Sn$--Frobenius algebra is essentially unique if it exists and
can be expressed solely in terms of the dual of the restriction maps
$\check r_{\tau}$ for $\tau$ a transposition. It thus depends only
on the maps $r_{\tau}$ and the metrics $\eta_{\tau}$. Also notice that
the $\Sn$ action acts transitively on all $A_{\tau}$ since any two
transpositions are conjugate.

If we specify the reconstruction data of the special $\Sn$ Frobenius manifold
to be that stemming from a symmetric product, it can furthermore
be shown that the unique structure does indeed exist \cite{K2}.

\subsection{Notation}
\label{sym}
 We define the degree of $\s \in \Sn$ to be
$|\s| :=$  the minimal length of $\s$ as a word in transpositions.
We define the length of $\s$ as $l(\s):=$ the number of cycles
in a cycle decomposition. Notice  $|\s|= n-l(\s)$.

We also set $l(\s_1, \dots,\s_n)= |\langle
\s_1,\dots,\s_n\rangle\backslash\{1,\dots,n\}|$ where $\langle
\s_1,\dots,\s_n\rangle$ is the group generated by the $s_i$ and
the quotient is by the natural permutation action.

\begin{df} We call two elements
$\sigma,\sigma' \in \Sn$ {transversal}, if
$|\sigma\sigma'|=|\sigma|+|\sigma'|$.
\end{df}
\subsection{Normalizability}
\label{norm}
\begin{df}
We call a non--abelian cocycle $\varphi$
normalized if
$\forall \t,\s\in \Sn, |\t|=1: \varphi_{\s,\t}= 1$.

We call a cocycle  $\g:\Sn \times \Sn \rightarrow A$
{normalizable}
if for all {transversal} pairs  $\t, \s \in \Sn , |\t|=1:
 \g_{\s,\t}\in A^*_e$, $A^*$ being the invertible elements of $A_e$, and
{normalized} if it is normalizable and
 for all {transversal } $\tau, \sigma \in \Sn , |\t|=1:
\g_{\s,\t}=1$.
\end{df}

In the example of symmetric products of an irreducible Frobenius
algebra the invertibles are precisely $k^*$.

\subsection{Discrete Torsion for the symmetric group}

As is well known (see e.g. \cite{Kar}) $H^2(\Sn,k^*)=\Zz$ and
$\Hom(\Sn,k^*)=\Zz$. We denote the non--trivial element of
$H^2(\Sn,k^*)$ by $[\Phi]$. There is a representative $\Phi$ of
this class which actually satisfies $\Phi(\s,\s')=\pm 1$ for
transversal $\s,\s'$. We denote the generator for the super--twist
by $\Sigma$.

\begin{thm}\cite{K2}
\label{normalize}
Any  non--abelian cocycle $\varphi$ after possibly twisting
by the discrete-torsion $\Phi$ and super-twist $\Sigma$ can be normalized.

Any compatible pair of a normalizable  $\Sn$ cocycle $\g$ and a
normalized non--abelian cocycle $\varphi$ can be normalized by a
rescaling $1_{\sigma} \mapsto \lambda_{\sigma}1_{\sigma}$.

And vice--versa given for
any normalized $\Sn$ cocycle $\g$ there are only
two compatible non--abelian cocycle $\varphi$
differing by the super--twist $\Sigma$ namely:

\begin{equation}
\varphi_{\s,\s'}= (-1)^{p|\s||\s'|}
\end{equation}
where $p$ is either $0$ or $1$.
\end{thm}

\subsection{Uniqueness}
\label{unique}
\begin{thm} \cite{K2}
 Given a special $\Sn$ algebra datum,
a choice of normalized cocycle $\g:\Sn \times \Sn \rightarrow A$
is unique up to a super--twist by $\Sigma$.

It is determined by $\g_{\t,\t}=\check r_{\t}(1_{\t})$ and without
the twist it is given by equation (\ref{expform}).
\end{thm}

\subsubsection{Explicit form of the cocycles}

\begin{equation}
\label{expform}
\g_{\s,\s'}=\pi_{{\s\s'}}(\g_{\s,\s'}\prod_{i=1}^{|\s'|}
\g_{\tau'_{i+1},\prod_{j=1}^{i}\tau'_{j}})
=\pi_{\s\s'}(\prod _{i=1}^{|\s'|}
\g_{\s\prod_{j=1}^{i}\tau'_{i-1},\tau'_{i}})
=\pi_{\s\s'}(\prod_{i \in I}\g_{\tau'_{i},\tau'_{i}})
\end{equation}
where
$I:=\{i: |\s(\prod_{j=1}^{i-1}\tau'_{j})\tau'_{i}|
=|\s\prod_{j=1}^{i-1}\tau'_{j}|-2\}$.

\subsection{Existence}

We would also would like to recall the following existence theorem:
\begin{thm} \cite{K2}
\label{existence}
The equations
\begin{eqnarray}
\label{nottrans}
r_{\s,\s'}(\g_{\s,\s'})&=&r_{\s\s'}(\prod_{i\in I}\g_{\t_{i},\t_{i}})=
\prod_{i\in I'}\pi_{\s\s'}(\g_{\t_{i},\t_{i}})\prod_{j\in I''}
r_{\s\s'}(\g_{\t_{j},\t_{j}})\nn\\
&=&\bar \g_{\s,\s'}\g_{\s,\s'}^{\perp}
\end{eqnarray}
where
\begin{eqnarray}
I'=\{ i \in I:
\pi_{\s,\s'}(\g_{\t_{i},\t_{i}})=\pi_{\s\s'}(\g_{\t_{i},\t_{i}})\}\nn\\
I''=\{ i \in I: \pi_{\s,\s'}(\g_{\t_{i},\t_{i}})
\neq\pi_{\s\s'}(\g_{\t_{i},\t_{i}})\}
\end{eqnarray}
and $\g^{\perp}_{\s,\s'}=\check r^{gh}_{g,h}(1_{g,h})$ are well
defined and yield a group cocycle compatible with the special
$\Sn$ intersection data.

$$A_{\s_1,\dots,\s_k}= (A^{\otimes l(\s_1,\dots,\s_k)},
\eta^{\otimes l(\s_1,\dots,\s_k)},
1^{\otimes l(\s_1,\dots,\s_k)})
$$
derived from any Frobenius algebra $(A,\eta,1)$. The restriction maps
are contractions via multiplication.
\end{thm}

For details on reconstruction data we refer to \cite{K1,K2}.
\subsubsection{Remark} In order to understand the form a the
cocycle above let us consider the case of the second symmetric
product. Thus we need to consider the ${\bf S}_2$-Frobenius
algebra:
$$A= A_e\oplus A_{\tau}= A\otimes A \oplus A$$
with the metric $\eta \otimes \eta \oplus \eta$. The restriction
map $r_{\tau}$ is just the multiplication $\mu: A \otimes A
\rightarrow A$. Now $\check r(1_{\tau})= \sum a_i \otimes b_i$ and
$r_{\tau} (\check r_{\tau})=\sum a_ib_i=e$ where $e$ is the Euler
class.

In general, if $\tau = (ij)$ then $r_{(ij)}$ contracts the i--th
and j--th component. Then $\g_{\t,\t}=\check r_{\tau}$ is a sum
over elements which differ from one only in the i--th and j--th
factor. Considering a product of $\g_{\t,\t}$ and restricting it
to $A_{\s\s'}$ amounts to performing several contractions using
the multiplication. For each individual $\g_{\t,\t}$ there are
only two choices: those which get contracted and yield Euler
classes -- this is the set $I''$-- and those that do not get
contracted -- this is the set $I'$.

\subsubsection{Remark} For the existence proof in \cite{K2} we used
the theory of intersection Frobenius algebras as can bee seen from
the decomposition of the cocycles. This also makes it easier
to compare with the results of \cite{LS}.

\subsubsection{Remarks}
\begin{itemize}
\item[1)] In the case that $A$ is the trivial one dimensional
Frobenius algebra this structure coincides with the group ring $k[\Sn]$.
\item[2)]
Applying our result to the situation where $A$ is the Frobenius
algebra associated to a variety or a compact space we recover
the results of \cite{FG} and taking invariants those of \cite{U}.
\end{itemize}

\begin{df} Given a Frobenius algebra $A$,
the series of $\Sn$--Frobenius algebras determined by Theorem
\ref{existence} is called the second quantization of $A$.
\end{df}

This terminology is based on \cite{DMVV}.

\section{The Twist for the Hilbert scheme}
If one changes the sign in the metric $\eta_{\tau}$
this also changes the multiplication. These changes are such, that
they are uniquely realized as a twist with a discrete torsion whose
cocycle $\a$ is determined by $\a(\tau,\tau)=-1$ and normalization.
This cocycle is actually trivial in cohomology with
coefficients in ${\bf C}^*$, but nevertheless changes the
multiplication and metric as desired.

\subsection{The twisted group ring $k^{\a}[\Sn]$}

Recall that $H^2(\Sn, {\mathbb C})= \ZZ$, but the twists are actually
given by $\a \in Z^2(\Sn,k)$.
Any normalized cocycle is fixed by $\a(\tau,\tau)$, $\tau$ any transposition.
Thus we may regard
$\a \in Z^2(\Sn,k)$ which is fixed by $\a(\tau,\tau)=-1$ for any
transposition $\tau$.

Notice that $[\a]=0 \in H^2(\Sn, {\mathbb C})$ as is easily seen from
the existence and uniqueness together with Remark \ref{conorm}.
It is  not trivial however in $H^2(\Sn, {\mathbb Q}^*)$.
In any case, we can consider the twist by this class.

\begin{prop}
\label{twist}
Given any $\Sn$--Frobenius algebra $A$, twisting it by the normalized
 discrete torsion $\a \in Z^2(G,k^*)$ defined by $\a(\t,\t)=-1$
 changes the structures via:

\begin{eqnarray}
a_{\s} \circ^{\a} b_{\s'} &=& (-1)^{\frac{1}{2}(|\s|+|\s|-|\s\s'|)}
a_{\s} \circ  b_{\s'}\nn\\
\varphi^{\a}_{\s}(a_{\s'}) &=&  \varphi_{\s}(a_{\s'})\nn\\
\chi_{\s}^{\a}&=&\chi_{\s}\nn\\
\eta^{\a}_{\s}&=&(-1)^{|\s|}\eta_{\s}
\label{metric}
\end{eqnarray}
\end{prop}

This is how the multiplication and metric get to be changed via a sign while
the $\Sn$ action remains unchanged.

\begin{proof} Due to \ref{defprop} it suffices to show that this holds for the
twisted group ring i.e.\ the following equations hold.
\begin{eqnarray}
\a(\s,\s') &=& (-1)^{\frac{1}{2}(|\s|+|\s|-|\s\s'|)}\\
\label{gamma}
\eps(\s,\s') &=&  1
\label{phi}
\end{eqnarray}
The equation (\ref{gamma}) follows
from $\a(\t,\t)=-1$ and the  general structure of $\a(\s,\s')$ of \ref{unique}
in particular from the equation (\ref{expform})
 by noticing that
$$
|I|=\{i: |\s(\prod_{j=1}^{i-1}\tau'_{j})\tau'_{i}|
=|\s\prod_{j=1}^{i-1}\tau'_{j}|-2\} =
\frac{1}{2}(|\s|+|\s'|-|\s\s'|)
$$

The correction of $\chi$ is always trivial and the one for $\varphi$
is given by $\eps$ which satisfies

$$
\eps(\s,\s') = \frac{\a(\s,\s')}{\a(\s\s'\s^{-1},\s)}
 = 1
$$
since
$$
|\s\s'\s^{-1}|=|\s'| \text { and so also } |\s'\s|= |\s\s'|
$$

Finally for the last equation of (\ref{metric}) we read off

$$\a(\s,\s^{-1})=
(-1)^{\frac{1}{2}(|\s|+|\s^{-1}|- |\s\s^{-1}|)}= (-1)^{|\s|}$$
\end{proof}

\subsection{A family of multiplications}
The only freedom of choice for a normalized cocycle is $\g_{\t,\t}$
which is determined uniquely from the metrics $\eta_e$ and $\eta_{\tau}$.

Given a fixed metric $\eta_{\t}$ it is  possible to change it by
homothety to $\eta^{\lambda}_{\t}= \lambda \eta_{\t}$ and keeping
the cocycles normalized using discrete torsion. This is achieved
by twisting with the normalized discrete torsion cocycle
determined by $\a(\t,\t)=\lambda$.

Vice--versa the only twists with discrete torsion cocycles that
keep the cocycle $\g$ normalized are fixed by their value
$\a(\t,\t)=\lambda$. The effect on the metric $\eta_{\tau,\tau}$
is a scaling by $\lambda$.

Using the same arguments as in \ref{twist}
\begin{prop}
Let $\a\in Z^2(\Sn,k^*)$ by the normalized cocycle determined by
$\a(\t,\t)=\lambda$  and $A$ by the $\Sn$--Frobenius algebra
associated to the symmetric product then the $\Sn$ Frobenius
algebra $A^{\a}$ is the twisted algebra found in \cite{QW}.
\end{prop}

This Proposition explains the existence and the special role of these
deformed multiplications. In the complex case, the existence
of this family also shows
the triviality of the cocycles $[\a]\in H^2(\Sn,{\bf C}^*)$.

\subsection{Remark} Notice that the discrete torsion $\eps$ is
trivial.
To obtain a non--trivial $\eps$ one could super--twist
\cite{K1,K2}, but this would also change the permutation action by
tensoring on the determinant representation, which is not intended
for the current application as we wish to keep symmetric invariants,
not anti--symmetric ones. Another way would be
to super--twist and to twist with a non--trivial cohomology class,
which would restore the
action to a symmetric one.
These results can be compared with \cite{D1,D2,DMVV,K2}.

\subsection{Application to the Hilbert scheme}

Comparing our analysis with that of \cite{LS} finally yields:

\begin{thm}{The cohomology of $Hilb^{[n]}$, the Hilbert scheme
of n--points for a K3 surface, is the $\Sn$ invariant part of
the unique $\Sn$--Frobenius Algebra
associated to the symmetric product of the cohomology of the surface
twisted by the specific discrete torsion given above. Or
in other words the $\Sn$--Frobenius Algebra structure for the extended
global orbifold cohomology twisted by specific discrete torsion, which is fixed
by the map of \cite{LS}.
In general, the sequence of spaces $Hilb^{[n]}$
gives rise to the twisted second quantization of the underlying surface
on the cohomological (motivic) level.}
\end{thm}

\subsubsection{Remark} The necessity of this type of twist and the particular
choice are fixed by two remarks.
First, the chosen isomorphism of \cite{LS} involves
the intersection form of the resolution of the double points,
which is negative definite. The intersection form of the natural
symmetric group Frobenius algebra
does not reflect this property.
The way that discrete torsion changes
the metric then fixes the cocycle as a normalized cocycle,
by regarding the metric of $A_{\tau}$ for $\tau$ a transposition.
This geometrically corresponds to a simple blow-up
along the diagonal.

\end{document}